\documentclass[12pt,a4paper]{article}
\usepackage{amsmath,amssymb,amsfonts}
\def\Bbb{\mathbb}

\title{\bf Cotangent power sums and character coordinates}

\author{Kurt Girstmair}

\date{}

\makeatletter
\let\@@maketitle=\maketitle
\def\maketitle{\def\thispagestyle##1{\relax}\@@maketitle}
\makeatother
\newtheorem{theorem}{Theorem}
\newtheorem{prop}{Proposition}

\def\BE{\begin{equation}}
\def\EE{\end{equation}}
\def\BD{\begin{displaymath}}
\def\ED{\end{displaymath}}
\def\BA{\begin{array}}
\def\EA{\end{array}}
\def\BEA{\begin{eqnarray*}}
\def\EEA{\end{eqnarray*}}
\def\BI{\bibitem}

\def\Z{\Bbb Z}
\def\Q{\Bbb Q}

\def\XX{{\cal X}}

\def\phi{\varphi}

\def\MB{\mbox}
\def\LD{\ldots}
\def\OV{\overline}

\def\WH{\widehat}

\def\sminus{\smallsetminus}
\def\DIV{\,|\,}

\def\BQ{``}
\def\EQ{'' }
\def\EQP{''}

\def\MN{\medskip\noindent}

\def\BRCHIF{B_{r,\chi_f}}

\def\BJCHIF{B_{j,\chi_f}}

\def\B3{B^{(3)}}
\def\ct3h{\WH{\rm{ct}}^{(3)}}

\newcommand{\btop}[2]{\genfrac{}{}{0pt}{1}{#1}{#2}}

\begin{document}
\maketitle

\begin{abstract}

\noindent
We show that certain sums studied in two recent papers are basically character coordinates (as they are called in the literature). These sums
involve values of Dirichlet characters and powers of $\cot(\pi k/n)$, $1\le k\le n-1$. We also show that a basic tool for the study of these sums
was already given in 1987, in the form of the character coordinates of so-called cotangent numbers. By means of this tool, we obtain the results of the said papers
in a simple and lucid way. We also show that the coefficients of the linear combinations used in the said papers are essentially the same.

\end{abstract}

\section*{1. Introduction}

Let $n\ge 2, r\ge 1$ be integers. Let $\zeta_n=e^{2\pi i/n}$ and $k\in\Z$, $(k,n)=1$. Since
\BE
\label{1.0}
i\cot(\pi k/n)=(1+\zeta_n^k)/(1-\zeta_n^k),
\EE
the numbers
$i^r\cot^r(\pi k/n)$  lie in the $n$th cyclotomic field  $\Q(\zeta_n)$.
Let $\chi$ be a Dirichlet character mod $n$. In two recent papers the character sums
\BE
\label{1.2}
  \sum_{k=1}^n\chi(k)i^r\cot^r(\pi k/n)
\EE
formed with these numbers have been expressed in terms of generalized Bernoulli numbers and  Gauss sums, see \cite[Cor. 13]{Is3}, \cite[Cor. 2.19]{Fr}.

However, the authors of these
papers do not refer to the fact that (\ref{1.2}) is basically a {\em character coordinate}, more precisely, the $\OV{\chi}$-coordinate of $i^r\cot^r(\pi/n)$.
Character coordinates have useful properties, see Sections 2 and 4. They have been known since 1959, see \cite{Leo} and Section 2.
In addition, they contain information about Galois modules, see \cite{Gi1}.

Moreover, the said authors make no use of the fact that the number $i^r\cot(\pi/n)$ has {\em natural components}  with respect to character coordinates.
namely, the so-called {\em cotangent numbers} $i^j\cot_{j-1}(\pi/n)$, $1\le j\le r$, where $\cot_l$ is the $l$th derivative of the function $\cot$ (in particular, $\cot_0=\cot$).
Indeed, for these components the character coordinates were given
by the present author already in 1987, see \cite{Gi1}. This fact has the following consequence. If we express
the function $i^r\cot^r$ as a rational linear combination of the functions $i^j\cot_{j-1}$, $1\le j\le r$, $j\equiv r$ mod $2$ (if $r$ is even, one must include the constant function $1$),
we immediately obtain the character coordinates of $i^r\cot(\pi/n)$, see Theorem \ref{t1}.
This theorem  is given in \cite{Is3}, whereas the paper \cite{Fr} has an equivalent result, but only for primitive characters (see Section 4).
We think that our approach to this theorem is the simplest one known so far.

Our plan is as follows. In Section 2 we recall some basic facts about character coordinates and exhibit the $\chi$-coordinates of the cotangent numbers.

In Section 3 we express  $i^r\cot^r$ as a linear combination of the functions $i^j\cot_{j-1}$, as described above (see Proposition \ref{p1}),
and obtain the said Theorem \ref{t1}.

In Section 4 we show that the result of \cite{Fr} gives the representation of $i^r\cot^r$ as a linear combination of the said functions in a different form,
see Theorem \ref{t2}.
Properties of character coordinates play a decisive role in this connection.

\section*{2. Character coordinates}

For a number $a\in\Q(\zeta_n)$ and a  Dirichlet character $\chi$ mod $n$,
the $\chi$-coordinate $y(\chi|a)$ is defined by
\BE
\label{2.2}
  y(\chi|a)\tau(\OV{\chi}_f)=\sum_{\btop{1\le k\le n}{(k,n)=1}}\OV{\chi}(k)\sigma_k(a),
\EE
see \cite{Leo}.
Here $f$ is the conductor of $\chi$, $\chi_f$ the character mod $f$ attached to $\chi$, $\OV{(\enspace)}$ the complex conjugation and
\BD
 \tau(\OV{\chi}_f)=  \sum_{k=1}^{f}\OV{\chi}_f(k)\zeta_{f}^{k}
\ED
the (primitive) Gauss sum;
furthermore, $\sigma_k$ is the Galois automorphism of $\Q(\zeta_n)$ defined by $\zeta_n\mapsto\zeta_n^k$, $(k,n)=1$.

Let $\Q(\chi)$ be the field of values of $\chi$. Then $y(\chi|a)\in \Q(\chi)$ and the map
\BD
 y(\chi|-):\Q(\zeta_n)\to \Q(\chi):a\mapsto y(\chi|a)
\ED is $\Q$-linear and $G$-invariant, so
\BD
  y(\chi|\sigma_k(a))=\chi(k)y(\chi|a).
\ED
From (\ref{1.0}) we obtain
\BE
\label{2.4}
  \sigma_k(i\cot(\pi/n))=i\cot(\pi k/n).
\EE
In view of (\ref{2.4}), the sum of (\ref{1.2})
has the form
\BD
 \sum_{\btop{1\le k\le n}{(k,n)=1}}\chi(k)i^r\cot^r(\pi k/n)=y(\OV{\chi}|i^r\cot^r(\pi/n))\tau(\chi_f).
\ED

A number $a\in\Q(\zeta_n)$ is uniquely determined by its character coordinates, as the reconstruction formula
\BD
   a=\frac 1{\phi(n)}\sum_{\chi\in\XX}y(\chi|a)\tau(\OV{\chi}_f)
\ED
shows, where $\XX$ is the set of all Dirichlet characters mod $n$ (see \cite{Leo}).

Let $\XX^+=\{\chi\in\XX;\chi(-1)=1\}$ and $\XX^-=\{\chi\in\XX;\chi(-1)=-1\}$.
If $a\in\Q(\zeta_n)$ is real, then $y(\chi|a)=0$ for all $\chi\in\XX^-$. If $a\in\Q(\zeta_n)$ is purely imaginary, then
$y(\chi|a)=0$ for all $\chi\in\XX^+$.

The cotangent number $i^r\cot_{r-1}(\pi/n)$ is real, if $r$ is even, and purely imaginary, if $r$ is odd. Accordingly,
$y(\chi|i^r\cot_{r-1}(\pi/n))$ vanishes in the cases \BQ $r$ is even, $\chi\in\XX^-$\EQ and \BQ $r$ is odd, $\chi\in\XX^+$\EQP.
In the remaining cases we have
\BE
\label{2.8}
  y(\chi|i^r\cot_{r-1}(\pi/n))=\frac{\chi(-1)(2n)^r}{rf^r}\prod_{p\DIV n}\left(1-\frac{\OV{\chi}_f(p)}{p^r}\right)\BRCHIF,
\EE
see \cite[Thm. 2]{Gi1}. Here
\BD
 \BRCHIF =f^{r-1}\sum_{k=1}^fB_r(k/f)\chi(k),
\ED
where $B_r(x)$ is the $r$th Bernoulli polynomial, see \cite[Prop. 4.1]{Wa}.

\section*{3. Cotangent powers and cotangent derivatives}

Let $k\ge 1$ be an integer.
Our main tool in this section is a special case of Lemma 4.1 in \cite{Is2} for certain functions in $t$, namely,
\BE
\label{3.2}
 \frac 1{(1-e^t)^k}=\frac 1{(k-1)!}\sum_{j=1}^k S(k,j)\frac{d^{j-1}}{dt^{j-1}}\left(\frac 1{1-e^t}\right),
\EE
where $S(k,j)$ is the absolute value of the respective Stirling number of the first kind (i.e., the number of permutations of $k$ objects with exactly $j$ cycles).
Since
\BD
  \frac 1{1-e^t}=\frac i2\cot(-it/2)+\frac 12,
\ED
we have, for $j\ge 2$,
\BD
  \frac{d^{j-1}}{dt^{j-1}}\left(\frac 1{1-e^t}\right)=(-1)^{j-1}(i/2)^j\cot_{j-1}(-it/2).
\ED
If we insert this into (\ref{3.2}), we obtain
\BE
\label{3.4}
 \frac 1{(1-e^t)^k}=\frac 12+\frac 1{(k-1)!}\sum_{j=1}^k S(k,j)(-1)^{j-1}(i/2)^j\cot_{j-1}(-it/2).
\EE
On the other hand,
\BD
 i\cot(-it/2)=\frac 2{1-e^t}-1.
\ED
Therefore, the binomial formula yields
\BD
 (i\cot(-it/2))^r=(-1)^r+\sum_{k=1}^r\binom{r}{k}\frac{2^k}{(1-e^t)^k}(-1)^{r-k}.
\ED
In this identity, we replace $1/(1-e^t)^k$ by the right hand side of (\ref{3.4}) and change the order of summation in the resulting double sum.
This gives
\BD
 (i\cot(-it/2))^r=C+\sum_{j=1}^ri^j\cot_{j-1}(-it/2)\frac{(-1)^{j-1}}{2^j}\sum_{k=j}^r\frac{(-1)^{r-k}2^k}{(k-1)!}\binom{r}{k}S(k,j),
\ED
with
\BD
 C=(-1)^r+\frac 12\sum_{k=1}^r \binom{r}{k}2^k(-1)^{r-k}=\frac{(-1)^r+1}2.
\ED
Thus, we may write
\BD
   (i\cot(-it/2))^r=\frac{(-1)^r+1}2+\sum_{j=1}^rc_{r,j}\, i^j\cot_{j-1}(-it/2),
\ED
with
\BE
\label{3.6}
  c_{r,j}=(-1)^{r-1}\sum_{k=j}^r\frac{(-2)^{k-j}}{(k-1)!}\binom{r}{k}S(k,j).
\EE
The change of variables $(-it/2)\mapsto x$ yields the following proposition.

\begin{prop} 
\label{p1}

For $r\ge 1$, we have

\BD
   i^r\cot^r=\frac{(-1)^r+1}2+\sum_{j=1}^rc_{r,j}\, i^j\cot_{j-1}
\ED
with $c_{r,j}$ as in \rm{(\ref{3.6})}.

\end{prop} 

Since the functions $i^j\cot_{j-1}$, $j\ge 1$, are $\Q$-linearly independent, we have $c_{r,j}=0$ if $r\not\equiv j$ mod $2$.
If $r$ is odd, we may write, therefore,
\BE
\label{3.8}
  i^r\cot^r=\sum_{\btop{1\le j\le r}{j\equiv 1\,\rm{mod}\, 2}}c_{r,j}\, i^j\cot_{j-1}.
\EE
If $r$ is even, we obtain,
\BE
\label{3.10}
  i^r\cot^r=1+\sum_{\btop{1\le j\le r}{j\equiv 0\,\rm{mod}\, 2}}c_{r,j}\, i^j\cot_{j-1}.
\EE

Let $\chi_0$ denote the principal character mod $n$.
We have, for $\chi\in\XX$,
\BE
\label{3.12}
 y(\chi|1)=\left\{
             \begin{array}{ll}
               0, & \hbox{if $\chi\ne\chi_0$;} \\
               \phi(n), & \hbox{if $\chi=\chi_0$.}
             \end{array}
           \right.
\EE

Since the $\chi$-coordinate is $\Q$-linear, the identity (\ref{2.8}) gives the following result.

\begin{theorem} 
\label{t1}

If $r\ge 1$ is odd and $\chi\in\XX^-$ has the conductor $f$, then
\BD
  y(\chi| i^r\cot^r(\pi/n))=-\sum_{\btop{1\le j\le r}{j\equiv 1\,\rm{mod}\, 2}}c_{r,j}\,\frac{(2n)^j}{jf^j}\prod_{p\DIV n}\left(1-\frac{\OV{\chi}_f(p)}{p^j}\right)\BJCHIF
\ED
with $c_{r,j}$ as in \rm{(\ref{3.6})}.

If $r\ge 2$ is even and $\chi\in\XX^+$ has the conductor $f$, then
\BD
  y(\chi| i^r\cot^r(\pi/n))=y(\chi|1)+\sum_{\btop{1\le j\le r}{j\equiv 0\,\rm{mod}\, 2}}c_{r,j}\,\frac{(2n)^j}{jf^j}\prod_{p\DIV n}\left(1-\frac{\OV{\chi}_f(p)}{p^j}\right)\BJCHIF
\ED
with $y(\chi|1)$ as in \rm{(\ref{3.12})} and $c_{r,j}$ as in \rm{(\ref{3.6})}.

\end{theorem} 

Theorem \ref{t1} is given in \cite[Cor. 13]{Is3} (for a preliminary version see \cite[Cor. 4.4]{Is}). We think, however, that our simple access to this theorem, i.e., 
via cotangent derivatives and the old result (\ref{2.8}), deserves to be noted.

\section*{4. Another form of the coefficients $c_{r,j}$}

In the paper \cite{Fr}, formulas for $y(\chi|i^r\cot^r(\pi/n))$, $\chi\in \XX$ primitive, are  given that involve coefficients $d_{r,j}$ seemingly different from the above $c_{r,j}$.
In this section we use properties of character coordinates in order to show that $c_{r,j}$ and $d_{r,j}$ are basically the same.

For an integer $r\ge 1$, we put
\BD
  \XX^r=\left\{
          \begin{array}{ll}
            \XX^-, & \hbox{if $r$ is odd;} \\
            \XX^+, & \hbox{if $r$ is even.}
          \end{array}
        \right.
\ED
Let $\chi\in\XX^r$ be a primitive character mod $n$. Formulas (2.29) and (2.30) of \cite{Fr} say, in our terminology,
\BE
\label{4.2}
  y(\chi|i^r\cot^r(\pi/n))=-2^r\sum_{\btop{1\le j\le r}{j\equiv r\,\rm{mod}\, 2}}d_{r,j}B_{j,\chi}/j!,
\EE
with
\BE
\label{4.3}
  d_{r,j}=\sum_{\btop{j_1,\LD,j_r\ge 0}{j+2j_1+\LD+2j_r=r}}\prod_{t=1}^r B_{2j_t}/(2j_t)!.
\EE
Here the numbers $B_{2j_t}$ are ordinary Bernoulli numbers.

On the other hand,
(\ref{2.8}) says, since $\chi$ is primitive,
\BD
\label{4.4}
  y(\chi|i^j\cot_{j-1}(\pi/n))=\frac{(-1)^r2^j}{j}B_{j,\chi}
\ED
for the numbers $j\equiv r$ mod $2$, $1\le j\le n$. Hence we may express the numbers $B_{j,\chi}$ of (\ref{4.2})
in terms of $y(\chi|i^j\cot_{j-1}(\pi/n))$. Thereby, we obtain
\BD
   y(\chi|i^r\cot^r(\pi/n))=\sum_{\btop{1\le j\le r}{j\equiv r\,\rm{mod}\, 2}}\frac{(-1)^{r+1}2^{r-j}}{(j-1)!}\,d_{r,j}\,y(\chi|i^j\cot_{j-1}(\pi/n))
\ED
The $\Q$-linearity of the $\chi$-coordinate yields
\BE
\label{4.6}
  y(\chi|i^r\cot^r(\pi/n))=y(\chi|\sum_{\btop{1\le j\le r}{j\equiv r\,\rm{mod}\, 2}}\frac{(-1)^{r+1}2^{r-j}}{(j-1)!}\,d_{r,j}\,i^j\cot_{j-1}(\pi/n))
\EE

Now suppose that $n=p$ is a prime.
First we assume that $r$ is odd.
Then all characters $\chi\in\XX^r$ are primitive. Thus, the identity
(\ref{4.6}) holds for all $\chi\in\XX^r$. However, for $\chi\in \XX^{r+1}$ both sides of (\ref{4.6}) vanish.
So (\ref{4.6}) is true for all $\chi\in\XX$.
In Section 2 we have seen that this means
\BD
  i^r\cot^{r}(\pi/p)=\sum_{\btop{1\le j\le r}{j\equiv r\,\rm{mod}\, 2}}\frac{(-1)^{r+1}2^{r-j}}{(j-1)!}d_{r,j}\,i^j\cot_{j-1}(\pi/p).
\ED
On the other hand, the identity (\ref{3.8}) implies
\BD
  i^r\cot^{r}(\pi/p)= \sum_{\btop{1\le j\le r}{j\equiv r\,\rm{mod}\, 2}}c_{r,j}i^j\cot_{j-1}(\pi/p).
\ED
We shall see below that, for a sufficiently large prime $p$, the numbers
$i^j\cot_{j-1}(\pi/p)$, $1\le j\le r$, $j\equiv r$ mod 2, are $\Q$-linearly independent.
Under this assumption, we may compare the coefficients on the right hand sides of the last two identities and get
\BE
\label{4.7}
  c_{r,j}=\frac{(-1)^{r+1}2^{r-j}}{(j-1)!}\,d_{r,j}, \enspace 1\le j\le r, j\equiv r \MB{ mod }2.
\EE

If $r$ is even, there are minor differences.  Again, let $n=p$ be a prime.
Then all characters $\chi\in\XX^r$ are primitive except the principal character $\chi_0$.
Hence (\ref{4.6}) holds  for all
$\chi\in\XX\sminus\{\chi_0\}$.
Suppose that
\BD
 y(\chi_0|\cot^{r}(\pi/p))=C_1\in\Q\: \MB{ and }\: y(\chi_0|\sum_{\btop{1\le j\le r}{j\equiv r\,\rm{mod}\, 2}}\frac{(-1)^{r+1}2^{r-j}}{(j-1)!}\,d_{r,j}\,i^j\cot_{j-1}(\pi/n))=C_2\in\Q.
\ED
Now we use (\ref{3.12}) and $\phi(p)=p-1$.
Therefore, the $\chi$-coordinate of $i^r\cot^{r}(\pi/p)-C_1/(p-1)$ agrees withe the $\chi$-coordinate of
\BD
  \sum_{\btop{1\le j\le r}{j\equiv r\,\rm{mod}\, 2}}\frac{(-1)^{r+1}2^{r-j}}{(j-1)!}\,d_{r,j}\,i^j\cot_{j-1}(\pi/n)-C_2/(p-1)
\ED
for each $\chi\in\XX$, and so these numbers are equal. Thus,
\BD
  i^r\cot^{r}(\pi/p)=\sum_{\btop{1\le j\le r}{j\equiv r\,\rm{mod}\, 2}}\frac{(-1)^{r+1}2^{r-j}}{(j-1)!}\,d_{r,j}\,i^j\cot_{j-1}(\pi/n)+(C_1-C_2)/(p-1).
\ED
In a similar way as above, we use the fact that the family
$(1;i^j\cot_{j-1}(\pi/p): 1\le j\le r, j\equiv r \MB{ mod }2)$ is $\Q$-linearly independent for a sufficiently large prime $p$.
The identity (\ref{3.10}) implies
\BD
  i^r\cot^r(\pi/n)=\sum_{\btop{1\le j\le r}{j\equiv r\,\rm{mod}\, 2}}c_{r,j}\, i^j\cot_{j-1}(\pi/n)+1.
\ED
On comparing the coefficients on the right hand side of the last two identities, we see that
(\ref{4.7}) holds also in this case.
Altogether, we have the following theorem.

\begin{theorem} 
\label{t2}

The coefficients $c_{r,j}$ and $d_{r,j}$ are connected by \rm{(\ref{4.7})}.

\end{theorem} 

We still have to show the aforesaid linear independence. By (\ref{2.4}),
the conjugates of $i\cot(\pi/p)$ are just the numbers
\BD
  i\cot(\pi k/p),\enspace 1\le k\le p-1.
\ED
Since the cotangent function is strictly monotonous in $(0,\pi)$, these numbers are pairwise different. Accordingly, the minimal polynomial of $i\cot(\pi/p)$ (over $\Q$) has the degree
$p-1$.
This means that the numbers
\BD
  i^j\cot(\pi/p)^j,\enspace j=0,\LD,p-2,
\ED
are $\Q$-linearly independent.

Now we choose $p$ such that $p-2\ge r$. If $r$ is odd, the families
\BD
 (i^j\cot^j: 1\le j\le r, j\equiv r \MB{ mod } 2) \MB{ and }(i^j\cot_{j-1}: 1\le j\le r, j\equiv r \MB{ mod } 2)
\ED
span the same $\Q$-vector space. This is also true for the families
\BD
 (i^j\cot^j(\pi/p): 1\le j\le r, j\equiv r \MB{ mod } 2) \MB{ and }(i^j\cot_{j-1}(\pi/p): 1\le j\le r, j\equiv r \MB{ mod } 2).
\ED
Accordingly, one of the latter families is $\Q$-linearly independent if, and only if, the other is $\Q$-linearly independent.
But we have shown the $\Q$-linear independence of the left hand family. So the right hand family is also $\Q$-linearly independent.
If $r$ is even, we work with the families
\BD
 (1;i^j\cot^j(\pi/p): 1\le j\le r, j\equiv r \MB{ mod } 2) \MB{ and }(1; i^j\cot_{j-1}(\pi/p): 1\le j\le r, j\equiv r \MB{ mod } 2)
\ED
instead.

\bigskip
\centerline{\bf Acknowledgment}

\MN
The author thanks Brad Isaacson for the important reference \cite{Is3}.

\bigskip
\centerline{\bf Conflicting interests and data availability}

\MN
The author declares that there are no conflicting interests. The manuscript has no associated data.


\MN
Kurt Girstmair\\
Institut f\"ur Mathematik \\
Universit\"at Innsbruck   \\
Technikerstr. 13/7        \\
A-6020 Innsbruck, Austria \\
Kurt.Girstmair@uibk.ac.at

\end{document}